\documentclass[10pt]{article}

\usepackage{amsmath}
\usepackage{amssymb}
\usepackage{amsthm}
\usepackage[matrix,curve,arrow]{xy}
\usepackage{amscd}
\usepackage{hyperref}

\newtheorem{thm}{Theorem}[section]
\newtheorem{lem}[thm]{Lemma}

\newtheorem{defin}[thm]{Definition}

\numberwithin{equation}{section}

 \newenvironment{note}[1][Note.]{\begin{trivlist}
    \item[\hskip \labelsep {\bfseries #1}]}{\end{trivlist}}

\newcommand{\fun}{\mathcal{C}^{\infty}(T^* \mathbb{R}^N)}
\newcommand{\rr}{\mathbb{R}}
\newcommand{\cc}{\mathbb{C}}
\newcommand{\nn}{\mathbb{N}}
\newcommand{\zz}{\mathbb{Z}}
\newcommand{\isol}{\bullet}
\newcommand{\Op}{\operatorname{Op}}
\newcommand{\ra}{\rightarrow}
\newcommand{\la}{\leftarrow}
\newcommand{\Lra}{\longrightarrow}
\newcommand{\Lla}{\longleftarrow}
\newcommand{\rra}{\rightrightarrows}
\newcommand{\da}{\downarrow}
\newcommand{\ua}{\uparrow}
\newcommand{\trace}{\operatorname{trace}}
\newcommand{\ka}{ \vdots k}
\newcommand{\wh}{\widehat}
\newcommand{\id}{\operatorname{id}}

\begin{document}
\title{The symbol of a function of an operator}
\author
{Alfonso Gracia--Saz \footnote{
Department of Mathematics.  University of California, Berkeley.
\newline \mbox{~~~~~} {\tt alfonso@math.berkeley.edu}
\newline \mbox{~~~~~}{\it 2000 MSC:} 53D55 (primary), 81S10 (secondary)
\newline \mbox{~~~~~}{\it Keywords:}  deformation quantization, Moyal product, symbol.}}
\maketitle

\begin{abstract}
We give an explicit formula for
the symbol of a function of an operator.  Given an   
operator $\widehat{A}$ on $L^2(\rr^N)$ with symbol $A \in \fun$ 
and a smooth function $f$, we obtain the symbol of $f(\widehat{A})$
 in terms of $A$.  As an application, Bohr--Sommerfeld quantization rules are explicitely calculated at order 4 in $\hbar$.
\end{abstract}

\section{Introduction.}
The goal of this paper is to give a realistically computable formula for 
the symbol of a function of an operator.  Let $\widehat{A}$ be an 
essentially 
self-adjoint operator in $L^2(\mathbb{R}^N)$
 with symbol (definitions to follow; see \textsection \ref{sec:bac})
$A \in
\fun$.  Let $f:\rr \rightarrow \rr$ be an analytic function.  Let 
$\widehat{B} 
= f(\widehat{A})$
 be an
operator with symbol $B$.  We want to write $B$ in terms of $A$.
We will derive the formula:
\begin{equation} \label{eq:primera}
B = \sum_{\Gamma} \left( \frac{i \hbar}{2} \right)^E 
	\frac{c_\Gamma}{S_\Gamma} \lambda_\Gamma(A) \;
	\frac{f^{(V)}(A)}{V!}
\end{equation}
	The sum is taken over finite graphs $\Gamma$ with no isolated
 vertices.  
For every such graph $\Gamma$,
 $V$ is the number of vertices and $E$ is the number of edges.
$\lambda_\Gamma(A)$ is a 
polynomial in the derivatives of $A$ constructed algorithmically from 
$\Gamma$ (see \textsection \ref{sec:gra}).  
$S_\Gamma$ is the order of 
the symmetry group of $\Gamma$.  $c_\Gamma$ is a simple invariant of 
$\Gamma$ (see \textsection \ref{sec:ng}).  The 
terms through order 4 in $\hbar$ of \eqref{eq:primera} are shown in 
Appendix \ref{sec:soafoaoao4ih}.

       The existence of a universal equation like \eqref{eq:primera} was used by Voros \cite{vor:aheosqs} and Colin de Verdi\`ere \cite{col:bsrtao} as part of a calculation to obtain Bohr--Sommerfeld quantization rules at higher orders in $\hbar$.  They derived it by using a explicit spectral theorem that writes $f(\wh{A})$ in terms of the resolvent.\footnote{In Appendix \ref{sec:aado} we give an alternative derivation of our main result inspired by this approach.}
Their method gives a recursive way of obtaining higher order corrections in $\hbar$ but is, in practice, intractable after order 2.  In contrast, the diagrammatic notation that we use (inspired by \cite{car:mspattbsa}) makes it simple to derive all our formulas and to write down explicitly their terms.

      For the derivation of \eqref{eq:primera} we need:
\begin{itemize} 
\item  Weyl quantization, i.e., a well defined correspondence between operators and symbols (see \textsection \ref{sec:wq}), but we will not use its explicit form.
\item  The explicit form of the Moyal product (see \textsection \ref{sec:tmp}).
\item  A spectral theorem, i.e., a way to define a function of an operator (see \textsection \ref{sec:st}), but we will not use its explicit form.
\end{itemize}

Before going on, let us mention some possible applications of this
calculation:
\begin{enumerate}
	\item   \emph{Bohr--Sommerfeld quantization rules.} This is treated in 
  \textsection \ref{sec:absr}.
        \item \emph{Determinant of certain differential operators.}  See \cite{loi:ymafmcboraot4dmp}.   Such determinants naturally arise in quantum field theory at the one loop level.  The determinant of an operator can be defined with the property $\log \det \wh{A} = \trace \log \wh{A}$, which holds in finite dimensional spaces.  Then, the trace of the operator $\wh{B} =\log \wh{A}$ can be calculated by means of its symbol $B$.
	\item \emph{Star exponential of quadratic forms}.  Omori et al
analyzed this problem in \cite{omo:sprtopiq}.  In the case when $A$ is a 
quadratic function and $f$ is an exponential, formula 
\eqref{eq:primera} simplifies.  See \textsection \ref{sec:tcoaqs}.
	 \end{enumerate}

       The structure of this paper is as follows.  In \textsection 
\ref{sec:bac} we explain the necessary background: Weyl quantization, the Moyal product and functions of operators.  In 
\textsection \ref{sec:gra} we introduce 
diagrammatic notation that will 
be used afterwards.  
     \textsection \ref{sec:cal} is the core of this paper, where  we derive \eqref{eq:primera} and other equivalent equations.
      In \textsection \ref{sec:gen} we consider generalizations to functions of various variables and to other quantizations.  
In \textsection \ref{sec:tcoaqs} we study the 
case of a quadratic symbol, in
particular the harmonic oscillator hamiltonian.  This is an example of restricting to a smaller class of symbols $A$, for which the family of graphs to consider becomes smaller, too, and \eqref{eq:primera} simplifies.
 \textsection \ref{sec:absr} explains the application to Bohr--Sommerfeld rules.

\begin{note}[Acknowledgments.]
	The 
author was supported by grants from La Caixa Foundation and
the Secretar\'ia de Estado de Universidades e Investigaci\'on del Ministerio Espa\~nol de Educaci\'on y Ciencia.  We thank the Institut de Math\'ematiques de Jussieu for their hospitality while part of this work was done.
  We are very thankful to Matthew Cargo, Laurent Charles, Yves Colin de Verdi\`ere, Robert Littlejohn and Alan Weinstein for useful, productive conversations.  We also thank David Farris, Alan Weinstein, Jared Weinstein and Marco Zambon for help with the editing of this paper.
\end{note}


\section{Background.} \label{sec:bac}

\subsection{Weyl quantization.} \label{sec:wq}

	Weyl quantization \cite{wey:guq} (or Weyl--Wigner correspondence) is a way to 
relate the classical and quantum descriptions of a system.   In the 
classical description of a system, the space of states is a Poisson 
manifold, whereas the quantum space is a Hilbert space.  For a particle 
in $N$ dimensions, the  Poisson 
manifold (classical) is $T^*\rr^N$ and the Hilbert space (quantum)
is $L^2(\rr^N)$.

	The observables are classically described by smooth functions on 
the phase space $\fun$.  In the quantum description they are operators on 
the Hilbert space $\Op(L^2(\rr^N))$.
For instance, we have canonical coordinates 
$$(x,p) := (x_1,\ldots,x_N,p_1,\ldots,p_N)$$ on the cotangent 
bundle $T^*\rr^N$: $x_j$ are the coordinates on the base 
$\rr^N$ and $p_j$ are the 
coordinates along the fibers.  $x_j$ and $p_j$ are elements of 
$\fun$ and hence 
classical observables.  We associate to them operators 
$\widehat{x_j}$ and $\widehat{p_j}$ on $L^2(\rr^N)$ defined by:
\begin{equation*}
\begin{split}
	\widehat{x_j} [\phi] (x) =& x_j \phi(x) \\
	\wh{p_j} [\phi] (x) =& -i \hbar \frac{d \phi}{dx_j}(x)
\end{split}
\end{equation*}
for $\phi \in L^2(\rr^N)$.

	Here $\hbar$ is the Planck constant.  We will treat it as a formal 
parameter and consider that all our spaces are formal power series in 
$\hbar$.  Thus the classical observables will be elements of 
$\fun[[\hbar]]$ and the quantum observables will be elements of 
$\Op(L^2(\rr^N))[[\hbar]]$.

        A \emph{quantization} is an extension of this correspondence to a 
map
\begin{equation*}
A \in \fun[[\hbar]] \; \mapsto \; \widehat{A} \in 
\Op(L^2(\rr^N))[[\hbar]]
\end{equation*}
  In particular the Weyl--Wigner correspondence is defined by:
\begin{equation}  \label{eq:weyl}
\begin{split}
\wh{A} [\phi](x) \; = \; & \int \frac{d^n y d^n p}{(2 \pi \hbar)^n}
e^{i(x-y)\cdot p/\hbar} A(\frac{x+y}2,p) \phi(y) \\  
A(x,p) \; = \; & \int \frac{d^n s}{(2 \pi \hbar)^n}
e^{-i s \cdot p / \hbar} \langle x+ s/2 \vert \wh{A} \vert
x-s/2 \rangle  
\end{split}
\end{equation}
  In order for it to be a well-defined  
bijection we need to 
restrict the domain and consider only a certain family of smooth functions 
and a certain family of operators. 
       See \cite{gri:mafdo} for details.  
In \textsection \ref{sec:oq} we consider alternative quantizations.

$A$ is called the \emph{symbol} of $\wh{A}$.  Moreover, the space of operators with composition is an algebra.  We can define an associative operation in $\fun[[\hbar]]$, called \emph{star product}  $\star$, that makes the bijection $A \mapsto \wh{A}$ into an algebra isomorphism.  In other words, $C \star D$ is the symbol of the operator $\wh{C} \wh{D}$.


\subsection{The Moyal product.} \label{sec:tmp}

Moyal \cite{moy:qmaast} gave an explicit expression for the star product in the case of the Weyl quantization, called the \emph{Moyal product} 
(but actually due to Groenewold \cite{gro:otpoeqm}).  It is derived from the definition of the Moyal product ($\wh{C\star D} = \wh{C}\wh{D}$) and the explicit form of Weyl quantization \eqref{eq:weyl}.
  If $C, D$ are symbols in $\fun$, then: 

\begin{equation} \label{eq:Moyal}
	C \star D = \sum_{k=0}^{\infty}
	\frac{1}{k!} \left( \frac{i \hbar}{2} \right)^k
	\{ C, D \}_k 
\end{equation}
	We need to explain the notation in \eqref{eq:Moyal}. 
In the natural 
coordinates $(x,p)$ on $T^*\rr^N$, 
the Poisson bivector field is
$$J = \sum_{j=1}^{N}
(\partial_{x_j} \otimes \partial_{p_j} - 
\partial_{p_j} \otimes \partial_{x_j})$$
where $\partial_q := \frac{\partial}{\partial q}$.  Let 
us call the coordinates 
$(z^1,\ldots,z^{2N}):=(x,p)$ and $J^{\mu \nu}$ the 
coefficients of the Poisson tensor on this chart 
\begin{equation*}
(J^{\mu \nu}) = \left( \begin{smallmatrix} 0 & I_N \\ -I_N & 0 
\end{smallmatrix} \right)
\end{equation*}
 where $I_N$ is the identity $N \times N$ matrix.
Using Einstein summation criterion (summation over repeated indexes),  
$J = J^{\mu \nu} \partial_\mu \otimes 
\partial_\nu$.  Let us also write
\begin{equation*}
  \begin{split}
C_{, \mu} & = \partial_{\mu} C 
= \partial_{z^\mu} C = \frac{\partial C}{\partial z^\mu} \\
C_{, \mu_1 \ldots \mu_k} & = \partial_{\mu_1} \ldots \partial_{\mu_k} C 
= \partial_{z^{\mu_1}} \ldots \partial_{z^{\mu_k}} C 
= \frac{\partial^k C}{\partial z^{\mu_1} \ldots \partial z^{\mu_k}}
  \end{split}
\end{equation*}
	Then the Poisson bracket in $T^*\rr^N$ can be written as
$$ \{C,D\} = C_{,\mu} J^{\mu \nu} C_{,\nu}$$
	The terms appearing in \eqref{eq:Moyal} are defined as:
\begin{equation*}
\begin{split}
	\{C,D\}_0 = & \;CD \\
	\{C,D\}_1 = & \; \{C,D\} = C_{,\mu} J^{\mu \nu} D_{,\nu} \\
	\{C,D\}_k = & \; C_{,\mu^1 \ldots \mu^k} J^{\mu^1 \nu^1} \cdots 
J^{\mu^k \nu^k} D_{,\nu^1 \ldots \nu^k}
\end{split}
\end{equation*}	

	Following \cite{car:mspattbsa} we use the notation 
``$\longrightarrow$'' for a Poisson tensor in the following way:   
``$\longrightarrow$'' is replaced by $J^{\mu \nu}$, the expression  
in the 
head of the arrow is acted on by $\partial_\nu$ and the expression in 
the 
tail of the arrow is acted on by $\partial_\mu$.  For instance
	\begin{equation*}
	\begin{split}
\xymatrix{C \ar[r] &D} = & \; C_{,\mu} J^{\mu \nu} D_{,\nu} \; = \; 
	\{C,D\} \\
\xymatrix{ C \ar@<.5ex>[r] \ar@<-.5ex>[r] & D} = & \;
	C_{,\mu^1 \mu^2} J^{\mu^1 \nu^1} J^{\mu^2 \nu^2} D_{,\nu^1 \nu^2}
	\; = \; \{C,D\}_2 \\
\xymatrix{ C \ar@<1.5ex>[r] \ar@<-1.5ex>^{\ka}[r] & D} = & \; 
	C_{,\mu^1 \ldots \mu^k} J^{\mu^1 \nu^1} \cdots
	J^{\mu^k \nu^k} D_{,\nu^1 \ldots \nu^k} \; =
	\; \{C,D\}_k
\end{split}
\end{equation*}
	Here, $\ka$ denotes k arrows.
	The same can be done with more complicated diagrams:
\begin{equation}  \label{eq:CDA}
\begin{split}
        \xymatrix{ C \ar[r] & D \ar@<-.5ex>[r] \ar[r] & A } 
        =& \quad \xymatrix{
        C \ar[r]^{J^{\mu_1 \nu_1}} & D
        \ar@<.5ex>[r]^{J^{\mu_2 \nu_2}}
	\ar[r]_{J^{\mu_3 \nu_3}} & A} \\
=& \quad   C_{,\mu_1} J^{\mu_1 \nu_1} D_{,\nu_1 \mu_2 \mu_3}
        J^{\mu_2 \nu_2} J^{\mu_3 \nu_3} A_{,\nu_2 \nu_3}
\end{split}
\end{equation}

	Since $J^{\mu \nu}$ is skew--symmetric, inverting an arrow multiplies 
the expression by $-1$:
\begin{equation*}
C \Lra D \; = \; (-1) \; C \longleftarrow D
\end{equation*}

	With this notation, the Moyal product \eqref{eq:Moyal} is written:
\begin{equation} \label{eq:moyala}
C \star D = \sum_{k=0}^{\infty}
        \frac{1}{k!} \left( \frac{i \hbar}{2} \right)^k
	\xymatrix{ C \ar@<1.5ex>[r] \ar@<-1.5ex>^{\ka}[r] & D}
\end{equation}

The fact that the Poisson bracket is a derivation on each argument 
(Leibniz rule or product rule) is written as:
\begin{displaymath}
\begin{array}{ccccc}
	(CD) \ra E &=& C \; (D\ra E) &+& (C \ra E) \;D \\
	\{ CD,E\} &=& C \{D,E\} &+& \{C,E\} D 
\end{array}
\end{displaymath}
Another application of the product rule is
\begin{equation} \label{eq:cadae}
\begin{array}{ccccc}
	(C \ra D) \ra E  &=& C \ra D \ra E &+& E \la C \ra D \\
	(C_{,\mu} J^{\mu \nu} D_{,\nu})_{,\alpha}
	J^{\alpha \beta} E_{,\beta} &=&
	C_{,\mu} J^{\mu \nu} D_{,\nu \alpha} J^{\alpha \beta} E_{,\beta}
	&+& 
	C_{,\mu \alpha} J^{\mu \nu} D_{,\nu} J^{\alpha \beta} E_{,\beta} 
\end{array}
\end{equation}
Notice that \eqref{eq:cadae} cannot be written with Poisson brackets 
$\{ \;,\;\} $.

	This notation makes it simpler to write certain calculations. 
 Lemma \ref{lem:diag1} in Appendix \ref{sec:alfcwg} is a generalization of 
\eqref{eq:cadae} that will be used in our derivations.

\subsection{Spectral theorems.} \label{sec:st}

      To understand \eqref{eq:primera} we need to define what a function of an operator means.  If $f(y)=y^n$ then $f(\wh{A}) = (\wh{A})^n$.  A spectral theorem extends this definition of  $f(\wh{A})$ to a wider class of functions $f$.  More specifically, let $\mathcal{A}$ be an algebra of smooth functions under pointwise multiplication and let $\wh{A}$ be an operator.  Then a \emph{spectral theorem} is an morphism of algebras 
$$ f \in \mathcal{A} \; \mapsto \; f(\wh{A}) \in \Op(L^2(\rr^N))$$
with certain properties.
     See \cite{dav:stado} for details.  There are many spectral theorems (i.e. different algorithms to calculate $f(\wh{A})$) for different algebras $\mathcal{A}$.  See \cite{and:fcfncowrsvaicf} for a complete list of references.
We will not use any explicit form of a spectral theorem (with the exception of
Appendix \ref{sec:aado}, where we derive an alternative proof of our main result). 
 We just need the fact that  a spectral theorem is a morphism of algebras:
$$f(\wh{A})g(\wh{A}) = (fg)(\wh{A})$$


\section{Graphs.} \label{sec:gra}

The formulas we are going to derive are power series whose terms are labeled by graphs.  We define now the family of graphs that we are going to use and introduce notation.

A \emph{graph} consists of a finite set of vertices and a finite number of edges.  Each vertex is represented by a dot. Each edge is represented by a line joining two vertices.  Multiple edges joining the same pair of vertices are allowed.  A self--edge (an edge from a vertex to itself) is not allowed.  A graph does not need to be connected.

An example of a graph is
\begin{equation} \label{eq:eins}
 \xymatrix{ \bullet \ar@{-}[r] \ar@/^/@{-}[rr] & \bullet \ar@{-}[r] &
\bullet \ar@{=}[r] & \bullet & \bullet \ar@{-}[r]& \bullet}
\end{equation}

A \emph{labeled graph} is a graph with in which we have labeled the vertices with the first $V$ natural numbers $1,2,\ldots, V$ and the edges with the first $E$ natural numbers $1,2,\ldots, E$.  Therefore, $V$ is the number of vertices and $E$ is the number of edges.  For instance:

$$ \xymatrix{ 1 \ar@{-}[r]_{1} \ar@/^/@{-}[rr]^{2} & 2 \ar@{-}[r]_{3} &
6 & 3 \ar@{-}@<-.5ex>[l]^{4} \ar@{-}[l]_{6} & 4 \ar@{-}^{5}[r] & 5}$$

  Even though a graph as defined above is not oriented, a labeled graph has a
natural orientation: every edge is oriented so that the target has a
higher label than the source:
$$ \xymatrix{ 1 \ar[r]_{1} \ar@/^/[rr]^{2} & 2 \ar[r]_{3} & 
6 & 3 \ar@<-.5ex>[l]^{4} \ar[l]_{6} & 4 \ar[r]^{5} & 5}$$

All the information in a graph is given by \emph{how many edges there are joining each pair of vertices}.
All the information in a labeled graph is given by \emph{which edge joins which pair of vertices}.  This leads us to adopt the following formal definitions.

\begin{defin}
Let $V$ and $E$ be two non--negative integers.  A \emph{labeled graph} with $V$ vertices and $E$ edges
is a map
$$s:  \{ 1, \ldots , E \} \Lra \mathcal{P}_2\{1, \ldots, V \} $$ 
where $\mathcal{P}_2X$ denotes the set of subsets of $X$ with 2 elements. 
\end{defin}
This means simply that the edge $i$ joins the pair of vertices $s(i)$.  The group $S_E$ of permutations of $E$ letters acts on $\{1,\ldots, E\}$.   The group $S_V$ of permutations of $V$ letters acts naturally on $\mathcal{P}_2\{1, \ldots, V\}$.  Hence the direct product $S_V \times S_E$ acts 
on functions $ \{ 1, \ldots , E \} \Lra \mathcal{P}_2\{1, \ldots, V \}$, that is,
on the set of labeled graphs with $V$ vertices and $E$ edges.
\begin{defin}
An \emph{unlabeled graph} or simply \emph{graph} with $V$ vertices and $E$ edges is an orbit of this action.
\end{defin}

        We will denote a labeled graph by $\Gamma$ and the corresponding (unlabeled) graph by $[\Gamma]$, if we need to distinguish between them.  Otherwise, we will abuse notation and denote a graph simply by $\Gamma$.

       It is convenient to define now two more concepts that will be needed later.  The \emph{order of symmetry} of a labeled graph is the number of permutations of edges and vertices that we can make without changing it.  Or, more formally:
\begin{defin}  \label{defin:S}
The \emph{order of symmetry} $S_\Gamma$ of a labeled graph $\Gamma$ is the order of the stabilizer of $\Gamma$ in the action of $S_V \times S_E$ on the set of graphs with $V$ vertices and $E$ edges.  
\end{defin}
  See Appendix \ref{sec:coc} for examples.

   A labeled graph is \emph{reduced} if it does not have any isolated vertices.  Or, more formally:
\begin{defin}  \label{defin:reduced} A labeled graph 
$s:  \{ 1, \ldots , E \} \Lra \mathcal{P}_2\{1, \ldots, V \}$
is \emph{reduced} if every $i=1,\ldots V$ is in some element of the image of $s$ (i.e., ``if every vertex is in some edge'').  
\end{defin}
Both concepts (order of symmetry of a graph and reduced graph) extend naturally to unlabeled graphs.
 All the previous examples are reduced.  The 
graph $\xymatrix{ \bullet \ar@{-}[r] & \bullet & \bullet}$ is not 
reduced. 

	Given a labeled graph $\Gamma$ with $V$ vertices, and given $V$ 
symbols $A_{1}, \ldots, A_{V}$ we construct a new symbol, called
$\lambda_\Gamma (A_1,\ldots,A_V)$,
 by substituting the vertices with $A_1, \ldots A_V$, 
and letting every 
edge represent a Poisson tensor (as explained in 
\textsection \ref{sec:tmp}).  We denote $\lambda_\Gamma (A, \ldots, A)$ 
simply by $\lambda_\Gamma (A)$.

	For instance, if $\Gamma$ is the labeled graph:
\begin{equation*}
\xymatrix{ 1 \ar@{-}[r]^{1} & 2 \ar@{=}[r]^{2}_{3}  & 3}
\end{equation*}
then $\lambda_\Gamma (C,D,A)$ is the expression in  \eqref{eq:CDA}.  
And
\begin{equation*}
\begin{split}
\lambda_\Gamma(A) \; = & \quad 
	\xymatrix{ A \ar[r] & A \ar@<-.5ex>[r] \ar[r] & A } \\
	=& \quad \xymatrix{
	A \ar[r]^{J^{\mu_1 \nu_1}} & A 
	\ar[r]^{J^{\mu_2 \nu_2}} \ar@<-.5ex>[r]_{J^{\mu_3 \nu_3}} & A} \\
=& \quad   A_{,\mu_1} J^{\mu_1 \nu_1} A_{,\nu_1 \mu_2 \mu_3}
	J^{\mu_2 \nu_2} J^{\mu_3 \nu_3} A_{,\nu_2 \nu_3}
\end{split}
\end{equation*}

	Since changing the direction of one arrow multiplies the expression by 
$-1$,  $\lambda_{[\Gamma]}$ is defined up to a sign.


\section{Main results and calculations.} \label{sec:cal} 
\label{sec:mrac}

We recall our problem. 
Let us fix an operator $\wh{A}$ with symbol $A$ and a smooth function $f$.  Let  $\wh{B} = f(\wh{A})$ be an operator with symbol $B$.  In this section we will perform the necessary calculations to obtain various expressions of $B$ in terms of $A$.

The main step is to obtain an expression for an iterated star product $C_1 \star \ldots \star C_n$ for symbols $C_i$ in terms of graphs.  We do this in \textsection \ref{sec:tnsp}.  Then in \textsection \ref{sec:ffftsoafoao} we derive our first expression for $B$ in terms of $A$.

Equation \eqref{eq:function} at the end of \textsection \ref{sec:ffftsoafoao} is a power series whose terms are parametrized by labeled graphs.  This is the easiest form of our result to derive, and it is useful for theoretical proofs.  However, it is not convenient for explicit calculations when we want to write the first few terms explicitely.  There are only a few unlabeled graphs, but many labeled graphs.  In \textsection \ref{sec:ng} we obtain our second expression for $B$ in terms of $A$, \eqref{eq:funred}, a series whose terms are parametrized by unlabeled graphs.  There is still a third form of our formula, Equation \eqref{eq:funexp}, whose terms are parametrized by connected graphs.  This last form is studied in \textsection \ref{sec:ef}.

Using either of these equations, we have included in Appendix \ref{sec:soafoaoao4ih}
the explicit form of the terms up to order $4$ in $\hbar$ of the symbol $B$ in terms of $A$.


\subsection{The $n$--th star product.} \label{sec:tnsp}

	The main step in the derivation of \eqref{eq:primera} is the 
following expression for the iterated star product, which generalizes Moyal's formula:

\begin{lem}  \label{lem:lemma}
	Let $C_1, \ldots, C_n \in \fun$ be symbols. Then
\begin{equation}	\label{eq:product}
  C_1 \star  \cdots \star C_n \;
= \; \sum_{k=0}^{\infty} \; \frac{1}{k!} 
\left(\frac{i \hbar}{2}\right)^k
	\sum_{\substack{\textrm{labeled graphs } 
	\Gamma \\ \textrm{ with } n \textrm{ vertices} \\
	\textrm{and $k$ edges} }} 
	\lambda_\Gamma (C_1, \ldots, C_n)
\end{equation}
\end{lem}
Note that when $n=2$, the previous lemma is exactly the Moyal formula for the star product \eqref{eq:moyala}.  For every $k$, there is only one labeled graph with $2$ vertices and $k$ edges.

   When $n=3$, for instance:
\begin{equation*}
\begin{split}
	A \star & A \star A \; = \quad A^3 \\
	& + \;  \frac{i \hbar}{2} \; \Big(
	 \xymatrix{A \ar[r] & A & A} \Big) 
	 \; + \; \frac{i \hbar}{2} \; \Big(
\xymatrix{A & A \ar[r]  & A} \Big) \\
	& + \; \frac{i \hbar}{2} \; \Big( 
	\xymatrix{A \ar@/^/[rr] & A & A} \Big) \; 
	+ \frac{1}{2!} \left( \frac{i \hbar}{2} \right)^2 \; \Big(
\xymatrix{A \ar[r] \ar@<-.5ex>[r] & A & A} \Big) \\
        &+ \;\frac{1}{2!} \left( \frac{i \hbar}{2} \right)^2 \; \Big(
        \xymatrix{A & A \ar[r] \ar@<-.5ex>[r]  & A} \Big) \; + \;
        \frac{1}{2!} \left( \frac{i \hbar}{2} \right)^2 \; \Big(
\xymatrix{A  \ar@/^/[rr] \ar@/^/@<.5ex>[rr] & A & A} \Big) \\
        &+ \; \frac{1}{2!}  
	\left( \frac{i \hbar}{2} \right)^2 \; \Big(
        \xymatrix{A \ar[r]^{1} & A  \ar[r]^{2} & A} \Big) \;  
	+ \; \frac{1}{2!} 
        \left( \frac{i \hbar}{2} \right)^2 \; \Big(
\xymatrix{A \ar[r]^{2} & A  \ar[r]^{1} &A} \Big) \\
	&+ \; \frac{1}{2!} 
        \left( \frac{i \hbar}{2} \right)^2 \; \Big(
	\xymatrix{A \ar[r]_{1} \ar@/^/[rr]^{2} & A & A} \Big) \;
        + \; \frac{1}{2!}
        \left( \frac{i \hbar}{2} \right)^2 \; \Big(
\xymatrix{A \ar[r]_{2} \ar@/^/[rr]^{1} & A & A} \Big) \\
	&+ \; \frac{1}{2!}  
	\left( \frac{i \hbar}{2} \right)^2 \; \Big(
        \xymatrix{A \ar@/^/[rr]^{1} & A \ar[r]_{2} & A} \Big) \; 
        + \; \frac{1}{2!}
        \left( \frac{i \hbar}{2} \right)^2 \; \Big(
\xymatrix{A \ar@/^/[rr]^{2} & A \ar[r]_{1} & A} \Big) \\
&+ \; O(\hbar^3) \\
	= \; & \quad A^3 \; + \; \frac{i \hbar}{2} 
	\left( \xymatrix{A \ar[r] & A } \right) A \;
	+ \; \frac{3}{2} \left( \frac{i \hbar}{2} \right)^2
\left( \xymatrix{A \ar[r] \ar@<-.5ex>[r] & A } \right) A\\
        &+ \; \frac{4-2}{2} \left( \frac{i \hbar}{2} \right)^2
	\Big( \xymatrix{A \ar[r] & A  & A \ar[l]} \Big) \;
+ \; O(\hbar^3)
\end{split}
\end{equation*}
     In the previous expression, the vertices are labeled from left 
to right $1,2,3$ in all graphs.

\begin{proof}[Proof of lemma \ref{lem:lemma}]
	We will use induction on $n$.  

	The result is true for $n=0,1,2$.

	Inductive step.  We use the associativity of the star product.
\begin{equation*}
\begin{split}
	C_1 \star \ldots& \star C_{n+1}  = 
(C_1 \star \ldots C_n) \star C_{n+1} \\
& = \left(
\sum_{\substack{\textrm{labeled} \\ \text{ graphs }
        \Gamma\\ \textrm{ with $n$  vertices}}}
	\; \frac{1}{E!} \left(\frac{i \hbar}{2}\right)^E
        \lambda_\Gamma (C_1,\ldots,C_n) \right) \star C_{n+1}		\\
& = \sum_{\substack{\textrm{labeled} \\ \text{ graphs }
        \Gamma \\ \textrm{ with $n$  vertices}}}
	\sum_{k=0}^{\infty}
        \; \frac{1}{E!k!} \left(\frac{i \hbar}{2}\right)^{E+k}
	\lambda_\Gamma(C_1,\ldots,C_n) \xymatrix{ \ar@<1.5ex>[r] 
	\ar@<-1.5ex>[r]^{\ka}&} C_{n+1}
\end{split}
\end{equation*}

	We can apply lemma \ref{lem:diag1} in Appendix \ref{sec:alfcwg}  to 
$\lambda_\Gamma(C_1,\ldots,C_n) \xymatrix{ \ar@<1.5ex>[r] 
\ar@<-1.5ex>[r]^{\ka}&} C_{n+1}$ and we will
get a sum over labeled graphs $\Gamma'$ with $n+1$ vertices.  They are 
built by starting with a labeled graph $\Gamma$ with $n$ vertices (labeled 
$1,\ldots,n$) and $E$ edges (labeled $1,\ldots,E$), 
adding the $(n+1)$-th vertex and $k$ edges 
(labeled $E+1,\ldots,E+k$) ending at the 
$(n+1)$-th vertex.  The number of edges of $\Gamma'$ is $E'=E+k$.  
      If we want to account for all possible labeled graphs with $n+1$ vertices and $E'$ edges, we need to divide by a factor of $\binom{E+k}{E}$ in order to account for all ways of relabeling the edges.
   Fortunately, we can write 
	$$\frac{1}{E!k!} = \frac{1}{(E+k)!} \binom{E+k}{E}$$
and we get
\begin{equation*}
        C_1 \star \ldots \star C_{n+1} \; =  
       \sum_{\substack{\textrm{labeled} \\ \text{ graphs }
        \Gamma' \\ \textrm{ with  $n+1$  vertices} }}
\frac{1}{E'!} \left(\frac{i \hbar}{2}\right)^{E'}
        \lambda_{\Gamma'}
\end{equation*}
\end{proof}


\subsection{First formula for the symbol of a function of an operator.} 
\label{sec:ffftsoafoao}

We now attack the problem of obtaining the symbol $B$ of $\wh{B}=f(\wh{A})$ in terms of $A$.  

To calculate $B$ we use the following expression \cite{cha:btoasca}:
\begin{lem} \label{lem:charles}
 Let $z_0 \in T^* \rr^N$ and let $a_0 := A(z_0)$.  Then
\begin{equation}  \label{eq:charles}
B(z_0) \; = \; \sum_{k=0}^{\infty} \;\frac{1}{k!} \; f^{(k)}(a_0) \; 
\left( A - a_0 \right)^{\star k} (z_0)
\end{equation}
\end{lem}
Note that the right hand side of \eqref{eq:charles} is well--defined only at the point $z_0$. 
To prove it, we need the following fact:
\begin{lem} \label{lem:nota}
Let $g$ be a smooth function and $\wh{C} := g(\wh{A})$.  Let $a_0 := A(z_0)$.
If $g$ has a zero of order $m$ at $a_0$, then $C(z_0)=O(\hbar^{m/2})$. \footnote{
Actually, we can do better: $C(z_0)=O(\hbar^{2m/3})$.
}
\end{lem}

\begin{proof}[Proof of lemma \ref{lem:nota}]
Let us write $g(y) = g_1(y) (y - a_0)^m$ and let $\wh{C_1} = g_1(\wh{A})$. 
Then $\wh{C} = g(\wh{A}) = g_1(\wh{A}) (\wh{A}-a_0)^m$ and 
$C = C_1 \star (A - a_0)^{*m}$.
Since $(A(z)-a_0)\vert_{z=z_0} = O(\hbar)$, then $C(z_0) = O(\hbar^{m/2})$.
\end{proof}

\begin{proof}[Proof of lemma \ref{lem:charles}]
For every $m$,  apply lemma \ref{lem:nota} with 
\begin{equation*}
g(y):= f(y) - \sum_{k=0}^{m} \frac{1}{k!} f^{(k)}(a_0) (y-a_0)^m 
\end{equation*}
  This proves that 
\begin{equation*}
B(z_0) \; - \; \sum_{k=0}^{\infty} \;\frac{1}{k!} \; f^{(k)}(a_0) \; 
\left( A - a_0 \right)^{\star k}(z_0) \; = O(\hbar^{m/2})
\end{equation*}
for all $m$.
\end{proof}

Now we only need to substitute \eqref{eq:product} into \eqref{eq:charles} to get an expression for $B$ in terms of graphs:
\begin{equation}
  B(z_0) \; = \; \sum_{n=0}^{\infty} \; \frac{1}{n!} \;
   f^{(n)}(A(z_0)) \; \sum_{k=0}^{\infty} \;
  \sum_{\substack{\textrm{labeled graphs } \Gamma \\ 
\textrm{with $n$ vertices}\\ \textrm{and $k$ edges}}}
\frac{1}{k!} \; \left( \frac{i \hbar}{2} \right)^{k}
\lambda_{\Gamma}(A-a_0) (z_0)
\end{equation}

In order to calculate $\lambda_{\Gamma}(A - a_0) (z_0)$ we need to write the symbol
$A-a_0$ at every vertex of $\Gamma$.  If a vertex is not isolated, then some derivatives are acting upon that symbol, and we may substitute $A-a_0$ with $A$, since $a_0$ is a constant.  If a vertex is isolated, then it 
contributes a factor of $(A-a_0)$,
and $(A-a_0)(z_0)=0$.  Hence we only need to consider graphs without isolated vertices, which we called reduced.  Putting this all together:

\begin{equation} \label{eq:function}
B = \sum_{\substack{\text{reduced labeled} \\ \text{graphs }
        \Gamma}} 
        \frac{1}{E!} \left(\frac{i \hbar}{2}\right)^E
	\frac{f^{(V)}(A)}{V!} \; \lambda_{\Gamma}(A)
\end{equation}

It is to be noted that the previous equation is not, strictly speaking, a power series expansion in $\hbar$, as $A$ itself depends on $\hbar$.  However, it is very simple to write $A = \sum_k \hbar^k A_k$ as a power series in $\hbar$ and expand \eqref{eq:function}.  As a matter of fact, we have an alternative way to write \eqref{eq:function}  with terms parametrized by \emph{graphs with weights} where every term is a monomial in $\hbar$.  We believe that the approach shown in this paper is simpler, though, and clearly illustrates the method.

Equation \eqref{eq:function} was easy to derive, and it is useful for proofs and theoretical calculations, as well as to generalize to other quantizations (see \textsection \ref{sec:oq}).  However, when we want to explicitly write the first few terms of it, this is not yet our ideal expression.  We can put together the contribution of labeled graphs that differ only in the labels to get a series whose terms are parametrized by unlabeled graphs.  We will do it next.

See Appendix \ref{sec:aado} for an alternative derivation of \eqref{eq:function}.


\subsection{Version with non--labeled graphs.}	\label{sec:ng}  \label{sec:vwngftc}

Two labeled graphs which are the same except 
for the labeling of vertices and 
edges give the same contribution (up to a sign).

\eqref{eq:function} can be rewritten as:

\begin{equation} \label{eq:cambio}
B = \sum_{\substack{\text{reduced}\\ \text{graphs }
        [\Gamma]}} 
        \frac{1}{E!} \left(\frac{i \hbar}{2}\right)^E
	\frac{f^{(V)}(A)}{V!}   
     \sum_{\Gamma' \in [\Gamma]}     
	\lambda_{\Gamma'}(A)
\end{equation}
    In words, we need to sum $\lambda_{\Gamma'}$ when $\Gamma'$ runs through 
all possible relabelings of $\Gamma$.  Remember that $\lambda_{\Gamma'}$ 
and $\lambda_{\Gamma}$ will be equal up to a sign.

Define
\begin{equation} \label{eq:defc}
	c_{\Gamma} = \sum_{\substack{\Gamma' \text{ is a reordering of} \\
	\text{the vertices of }	\Gamma}} 
	(-1)^{\# \text{ of arrows inverted going from $\Gamma$ to $\Gamma'$} }
\end{equation}
That is, we start with a labeled graph $\Gamma$.  Then, we consider the $V!$ possible ways of numbering the vertices of 
the graph with
$1,2,\cdots, V$.  For each of them, we orient the arrows so that 
they all go from the vertex with the lowest label to the vertex with the 
highest label.  Then we count these $V!$ relabelings with a sign, 
depending on the parity of the  number of arrows inverted from our original orientation.

	We can then write the contribution from \eqref{eq:cambio} 
as:
\begin{equation}  \label{eq:sim}
	\sum_{\Gamma' \in [\Gamma]} \lambda_{\Gamma'}(A) =
	\frac{E!}{S_{\Gamma}} \; c_{\Gamma} \; \lambda_{\Gamma}(A)
\end{equation}
	where $S_{\Gamma}$ is the order of the symmetry group of the 
(unlabeled) graph $[\Gamma]$ (see Definition \ref{defin:S}).  
The contribution corresponding to 
different relabelings of the edges is $E!$.  The contribution 
corresponding to different relabelings of the vertices is in $c_\Gamma$.  
And 
we have to divide by the order of the symmetry group, to
account for the situation in 
which exchanging edges or vertices results in the same labeled graph.

	For instance, if $\Gamma$ is the graph  $\bullet \ra \bullet \ra 
\bullet$, then the contribution from renumbering the vertices is:
\begin{displaymath}
\begin{array}{cccccccc}
	&& 1 \ra 2 \ra 3 && 1 \ra 3 \la 2 && 2\la 1 \ra 3 \\
	c_\Gamma &= & (-1)^0 & + & (-1)^1 & + & (-1)^1 \\ \\
	&& 2 \ra 3 \la 1 && 3 \la 1 \ra 2 && 3 \la 2 \la 1 \\
	& + & (-1)^1 & + & (-1)^1 & + & (-1)^2 & = \; -2
\end{array}
\end{displaymath}

	Finally we just have to substitute \eqref{eq:sim} into
 \eqref{eq:cambio}:
\begin{equation} \label{eq:funred}
B = \sum_{\substack{\text{reduced}\\ \text{graphs }
        [\Gamma]}}
        \left(\frac{i \hbar}{2}\right)^E
        \frac{c_\Gamma}{S_\Gamma}
	\frac{f^{(V)}(A)}{V!}
	\lambda_\Gamma(A)
\end{equation}

	Notice how $c_\Gamma$ and $\lambda_\Gamma(A)$ are only defined up to
a sign for $[\Gamma]$.  However, those signs cancel in their product $c_\Gamma
\lambda_\Gamma(A)$, which is well defined.

The explicit calculation of $c_\Gamma$ for a particular graph is actually 
very easy.  See Appendix \ref{sec:c},
which includes the value of $c_\Gamma$ and $S_\Gamma$ 
for all 
reduced, connected graphs with 2 or 4 edges.  Thanks to lemma \ref{lem:odd} 
we only need to consider graphs 
$\Gamma$ where every connected component has an even number of edges, since
otherwise $c_\Gamma = 0$.


\subsection{Version with connected graphs.} 
   \label{sec:ef} \label{sec:efftc} \label{sec:vwcg}

  Let us rewrite \eqref{eq:funred} as:
\begin{equation} \label{eq:D}
B = \left[ \sum_{\substack{\text{reduced}\\ \text{graphs }
        [\Gamma]}}
        \left(\frac{i \hbar}{2}\right)^E
	\frac{c_\Gamma}{ S_\Gamma} \;
	\lambda_\Gamma(A) \;
	\frac{D^V}{V!} \right] \; f(A)
\end{equation}
where $D$ is the differential operator which applies to $f$.

Whenever we have an expression like \eqref{eq:D}, a series labeled by a 
certain family of diagrams, it is standard to 
reduce all calculations 
to only connected diagrams.   Let's generalize.  

Let $\mathcal{G}$ be the free commutative monoid generated by the set 
$\mathcal{G}_0$.  Let $S$ be a commutative ring (with 
multiplicative notation).  Let 
$\mathcal{O}:\mathcal{G} \rightarrow S$ be a map satisfying
\begin{equation} \label{eq:propexp}
	\mathcal{O}(r_1 x_1 + \cdots + r_n x_n) \; = \;
	\frac{1}{r_1!\cdots r_n!}
	\left( \mathcal{O} (x_1)\right)^{r_1} \cdots
	\left( \mathcal{O} (x_n)\right)^{r_n}
\end{equation}
for distinct $x_1, \ldots, x_n \in \mathcal{G}$ and 
$r_1,\ldots,r_n \in \nn$.
Then, formally:
\begin{equation} \label{eq:expsum}
\sum_{x \in \mathcal{G}} \mathcal{O}(x) =
	\exp \left[ 
  \sum_{x \in \mathcal{G}_0} \mathcal{O}(x)
		\right]
\end{equation}

	In particular, consider $\mathcal{G}$ to be a family of diagrams 
closed 
under topological sum and generated by the connected non-empty diagrams 
$\mathcal{G}_0$.  If we write $\mathcal{O}(\Gamma) = 
\frac{\mathcal{M}(\Gamma)}{S_\Gamma}$, where $S_\Gamma$ is the order of 
the symmetry group of the diagram and $\mathcal{M}$ is a multiplicative 
function
\begin{equation*}
	\mathcal{M}(x_1 + x_2) =
	\mathcal{M}(x_1) \mathcal{M}(x_2)
	\text{  for all  } x_1, x_2 \in \mathcal{G}
\end{equation*}
then  $\mathcal{O}$ satisfies \eqref{eq:propexp}.

	In our case, 
\begin{equation*}
\begin{split}
& \mathcal{G} = \text{ reduced graphs} \\
& \mathcal{G}_0 = \text{ reduced, connected, non-empty graphs} \\
& S = \fun \\
& \mathcal{M} ([\Gamma]) =         \left(\frac{i \hbar}{2}\right)^E
        {c_\Gamma} \; \lambda_\Gamma(A) \;
        \frac{D^V}{V!}
\end{split}
\end{equation*}
$\mathcal{M}$ is multiplicative from \eqref{eq:cdis} in Appendix \ref{sec:c}
\begin{equation*}
	\frac{c_{\Gamma}}{V!} \; = \; 
	 \frac{c_{\Gamma_1}}{V_1!} \cdots 
	\frac{c_{\Gamma_n}}{V_n!}
\end{equation*}
	Hence, using \eqref{eq:expsum} in \eqref{eq:D}:
\begin{equation} \label{eq:funexp}
B = \left[ \exp 
	\sum_{\substack{\text{connected, reduced,}\\ \text{non-empty 
	graphs }[\Gamma]}}
        \left(\frac{i \hbar}{2}\right)^E
        \frac{c_\Gamma}{S_\Gamma} \;
        \lambda_\Gamma(A) \;
        \frac{D^V}{V!} \right]  \; f(A)  
\end{equation}


\section{Generalizations.} \label{sec:gen}

As we mentioned in the introduction, we needed three things to derive our results:
\begin{itemize}
\item Weyl quantization.
\item The Moyal product.
\item A spectral theorem.
\end{itemize}
Of these, we only used the explicit form of the Moyal product.  Actually, the form of the Moyal product is calculated from the form of the Weyl quantization.  We discuss now how to generalize to functions of various variables and other quantizations.


\subsection{Functions of several variables.} \label{sec:fovv} \label{sec:fonv}
\label{sec:fosv}

Let $\wh{A_1},\dots,\wh{A_n}$ be $n$ commuting operators in $L^2(\rr^N)$ with symbols $A_1,\ldots,A_n$.
Let $F:\rr^n \Lra \rr$ be a smooth function.  We consider the operator 
$\wh{B} = F(\wh{A_1},\ldots,\wh{A_n})$ with symbol $B$.  Can we extend our results to calculate $B$ in terms of $A_1,\ldots,A_n$?  The answer is yes.  We are still using Weyl quantization and the Moyal product, and spectral theorems  behave equally well for functions with several variables.  See \cite{and:fcfncowrsvaicf} for a list of references.

Hence, we only need to repeat our calculations, starting from \textsection \ref{sec:ffftsoafoao}, but with a function of a several variables.
 The counterpart of \eqref{eq:funred} is
    \begin{equation} \label{eq:n} 
B = \sum_{\substack{\text{reduced}\\ \text{graphs }
        [\Gamma]}}
        \left(\frac{i \hbar}{2}\right)^E
        \frac{c_\Gamma}{S_\Gamma} \;
        \frac{\partial_{i_1} \cdots \partial_{i_V} F(A)}{V!} \;
        \lambda_\Gamma (A_{i_1}, \cdots A_{i_V})
\end{equation}

For instance the first few terms of \eqref{eq:funred} are
\begin{equation*}
	B = f(A) - \frac{\hbar^2}{4} 
	\left[
	\frac{1}{2} \xymatrix{A \ar[r] \ar@<-.5ex>[r]&A}
		\frac{f''(A)}{2!} \quad + \quad
	A \ra A \la A \; \frac{f'''(A)}{3!}
	\right] + O(\hbar^4)
\end{equation*}
	and the first few terms of \eqref{eq:n} are
\begin{equation*}
\begin{split}
	B \;=\;& F(A) \\
	 & - \frac{\hbar^2}{4}
        \left[
        \frac{1}{2} \xymatrix{A_i \ar[r] \ar@<-.5ex>[r]&A_j}
                \frac{\partial_i \partial_j F(A)}{2!} \quad + \quad
        A_i \ra A_j \la A_k \; 
	\frac{\partial_i \partial_j \partial_k F(A)}{3!}
        \right] \\ &+ O(\hbar^4)
\end{split}
\end{equation*}

There are also the obvious versions with labeled or connected graphs.


\subsection{Other quantizations.} \label{sec:oq}

There are other quantizations apart from Weyl quantization, that is, correspondances between operators and symbols \cite{hir:dqittoqm}.  If we want to use them, then we have a different star product instead of the Moyal product.  The explicit form of the Moyal product has been used in two places: to prove lemma \ref{lem:nota} and to derive an expression for the iterated star product (lemma \ref{lem:lemma}).

Let us consider a generic star product that has the form 
\begin{equation} \label{eq:star}
	C \star D = \sum_{k=0}^{\infty}
	\frac{1}{k!} \left( \frac{i \hbar}{2} \right)^k
	\{ C, D \}_k 
\end{equation}
As long as $(C,D) \mapsto \{C,D\}_k$ is a bidifferential operator of order $m_k$ and $\lim_{k \ra \infty} m_k = \infty$, then lemma \ref{lem:nota} is satisfied.  Hence, for those star products, we only need to obtain an analogue to lemma \ref{lem:lemma}, that is, an expression for the iterated star product $C_1 \star \ldots \star C_n$ in terms of diagrams.  This can often be done by induction if we start by writing the star product of two symbols as a series in terms of diagrams.

For instance, we can consider standard order quantization:
\begin{equation*}
   A(x,p) = u_{-p} (x) \wh{A}[u_p](x)
\end{equation*}
where $u_p (x) = e^{ip \cdot x/\hbar}$.  In that case, the star product has the form of \eqref{eq:star} with $\{ C , D \}_k := \partial_{p_j}^k C \partial_{x^j}^k D$.  Lemma \ref{lem:nota} still holds and, if we change the meaning of ``$\Lra$'', lemma \ref{lem:lemma} is also true.  The simplest way is to use the same definitions we gave in \textsection \ref{sec:tmp} and \textsection \ref{sec:gra} to construct a polynomial $\lambda_\Gamma(A)$ from a labeled graph 
$\Gamma$, but using the (non--Poisson) tensor 
\begin{equation*} (J^{\mu \nu}) = \left( 
\begin{smallmatrix} 0 & 0 \\ I_N & 0  \end{smallmatrix} \right)
 \end{equation*}
With that convention, \eqref{eq:function} is still valid.  However, 
it is no longer true that $C \Lra D = - \; C \Lla D$.  As a consequence,
the versions with  unlabeled or connected graphs are messier (although they still exist).

In \cite{kon:dqpm}, Kontsevich gave a star product that quantizes any Poisson structure on $\rr^N$.  His expression is already a power series in $\hbar$ whose terms are labeled by a family of diagrams.  We can use it to obtain an equivalent of lemma \ref{lem:lemma} and we can derive, again, the counterparts of Equations \eqref{eq:function}, \eqref{eq:funred} and \eqref{eq:funexp}.  In fact, our set of labeled graphs is a subset of Kontsevich's set of labeled graphs.  It is to be noted, though, that Kontsevich's star product includes a weight $\omega_\Gamma$ associated to every diagram $\Gamma$ which is, in practice, hard to calculate. 
(See \cite{kat:kuffdqatcbhf} and \cite{pol:qolpsadom} for some results.) 


\section{The case of a quadratic symbol.} \label{sec:tcoaqs}

When we restrict to a smaller class of symbols, it is possible that the contribution of many graphs vanishes, simplifying our calculations.  As an example, we study here quadratic symbols.

	A particular (simple) case of importance consists of taking 
$\widehat{A} = \widehat{I}$, the harmonic oscillator hamiltonian in 1 dimension.
In the standard coordinates $(z^1,z^2)=(x,p)$ its symbol is 
$I=\frac{1}{2}((z^1)^2+(z^2)^2)$.  This simplifies the calculations 
because any 
third derivative vanishes:
$I_{,\mu_1\mu_2 \mu_3} = 0$.  The same is true for any quadratic 
function.   Let $Q_{\mu \nu}$ be a $2\times 2$ real symmetric matrix and 
consider the symbol $A = \frac{1}{2} z^{\mu} Q_{\mu \nu} z^{\nu}$.  Assume 
$A$ is the symbol of an operator $\widehat{A}$.  Then:
	\begin{equation*}
\begin{split}
	A_{,\nu} & = z^{\mu} Q_{\mu \nu} \\
	A_{,\mu \nu} & = Q_{\mu \nu}
\end{split}
\end{equation*}
and any third derivative vanishes.
	As a consequence, we only 
need to consider graphs where every vertex has at most two edges. After 
 lemma \ref{lem:odd} in Appendix \ref{sec:c}, 
we only need to consider graphs with an even number of edges.
 If we 
also ask them to be connected, reduced and non-empty, there are only two 
families of such graphs:
\begin{displaymath}
\begin{array}{lcccccccccc}
	 \Delta_k  &\; =\;&  \bullet &\ra& \bullet& \ra&
	 \cdots &\ra& \bullet&  \text{  (2k edges)}& k \geq 1 \\
	 \Lambda_k  &\; =\;& \bullet &\ra& \bullet& \ra&
	 \cdots &\ra& \bullet&  \text{  (2k edges)}& k \geq 1 \\
	&& \ua &&&&&&\da \\
	&&\bullet &\la& \bullet& \la& \cdots &\la& \bullet 		
\end{array}
\end{displaymath}

        In words, $\Delta_k$ consists of $2k+1$ vertices and $2k$ 
edges joined forming a line.  Its symmetry group has order $2$.
 $\Lambda_k$ consists of $2k$ vertices and $2k$ edges joined forming a 
simple cycle.  Its symmetry group has order $4k$.

	The corresponding polynomials in the derivatives of $A$ can be 
calculated:  
\begin{equation*}
\begin{split}
	\lambda_{\Delta_k} = &\; A_{,\mu_1} J^{\mu_1 \nu_1} 
A_{,\nu_1 \mu_2}
J^{\mu_2 \nu_2}  \cdots J^{\mu_{2k} \nu_{2k}} A_{,\nu_{2k}} \\
	= &\; z^{\mu_1} ((QJ)^{2k}Q)_{\mu_1 \nu_{2k}} z^{\mu_{2k}} \\
	\lambda_{\Lambda_k} = &\; A_{,\nu_{2k} \mu_1} J^{\mu_1 \nu_1} 
A_{,\nu_1 \mu_2} J^{\mu_2 \nu_2}  \cdots J^{\mu_{2k} \nu_{2k}} = \\
	= &\;\trace \;((QJ)^{2k}) 
\end{split}
\end{equation*}
	Since $Q$ is symmetric and $J$ is skew--symmetric, their product 
$QJ$ is traceless.  Hence $(QJ)^2 = -\det (QJ) \id$.  Write 
$ \omega^2 := \det (QJ) = \det(Q)$.  Then:
\begin{equation*}
\begin{split}
	\lambda_{\Delta_k} = & (-1)^k \omega^{2k} z^{\mu} Q_{\mu \nu} 
z^{\nu} = (-1)^k \omega^{2k} 2 A \\
	\lambda_{\Lambda_k} = & (-1)^k \omega^{2k} \trace({\id}) =
(-1)^k \omega^{2k} 2
\end{split}
\end{equation*}

	As for the coefficients $c_{\Delta_k}$ and $c_{\Lambda_k}$, Fact 
\ref{it:p} in Appendix \ref{sec:c} gives us the relation 
$c_{\Lambda_k} = -2k c_{\Delta_{k-1}}$.

	We now plug all this values in \eqref{eq:funexp} to obtain, 
for a general function $f$ and a quadratic symbol $A$:
\begin{equation} \label{eq:Inc}
\begin{split}
B \; = \Bigg(  \exp \Bigg[
	A & \sum_{k=0}^{\infty}
	\left( \frac{i \hbar}{2} \omega \right)^{2k} 
	\; \vert c_{\Delta_{k+1}} \vert \;
	\frac{D^{2k+1}}{(2k+1)!} \\
 & + \;	 \sum_{k=0}^{\infty}
        \left( \frac{i \hbar}{2} \omega \right)^{2k} 
	\; \vert c_{\Delta_k} \vert \;
	\frac{D^{2k}}{(2k)!}
\Bigg] \; \Bigg) \quad f(A)
\end{split}
\end{equation}

	We are left with a combinatorics problem: the sequence 
$\{ c_{\Delta_k} \}$. 
  Fact \ref{it:p} in 
Appendix \ref{sec:c} gives us a recurrence formula:
\begin{equation}	\label{eq:crec}
	c_{\Delta_k} = - \sum_{j=0}^{k-1} \binom{2k}{2j+1} \; 
c_{\Delta_j} \; c_{\Delta_{k-j-1}}
\end{equation}
The sequence is alternating in sign, and the first few 
absolute values are 1, 2, 16, 272, 7936 \ldots  This sequence is called 
the \emph{Zag numbers} \cite{slo:toeois}
and they appear in the McLaurin expansion of the tangent:
\begin{equation} \label{eq:Zag}
	\tan x = \sum_{k=0}^{\infty} 
\frac{|c_{\Delta_k}|}{(2k+1)!} x^{2k+1} 
\end{equation}

	To prove this, notice that $\tan x$ is the only odd solution to 
the differential equation 
\begin{equation}	\label{eq:de}
y' = 1 + y^2
\end{equation}	
	  Write a generic solution of the 
form $y = f(x) = \sum_{k=0}^{\infty} \frac{\alpha_k}{(2k+1)!} x^{2k+1}$ 
and substitute it into \eqref{eq:de}.  Equating coefficients, we conclude 
that the sequence 
$\alpha_k$ satisfies the same recurrence relation as the sequence 
$c_{\Delta_k}$ (Equation \eqref{eq:crec}).  Hence $\alpha_k = 
\vert c_{\Delta_k} \vert$.

	The Zag numbers can be written in terms of the Bernoulli numbers 
$B_n$:
\begin{equation*}
	|c_{\Delta_k}| = \frac{2^{2k}}{2^{2k}-1} \frac{|B_{2k}|}{2k}
\end{equation*}

	When we use equation \eqref{eq:Zag} in \eqref{eq:Inc}
 we obtain a nice, compact expression:
\begin{equation} \label{eq:Ic}
\begin{split}
	\widehat{B} \;= & \; f(\widehat{A}) \\
	B \; = \;& \sec \frac{i \hbar \omega D}{2} \; \exp \left[ 
	\frac{2A}{i \hbar \omega } \tan \frac{i \hbar \omega D}{2} -AD
\right] \; f(A)
\end{split}
\end{equation}

	Remember that $D$ is the derivative operator that applies to $f$.  
 If we also take the function 
$f(y) = e^{\varepsilon y}$, then $D$ acts simply as multiplication by 
$\varepsilon$.  In particular, when we consider the time evolution operator:
\begin{equation} \label{eq:bufalos}
\begin{split}
	\widehat{B} \; = & \; e^{-i t \wh{A} / \hbar} \\
	B \; = & \; \sec \frac{t \omega}{2} \; \exp 
	\left[ 	 \frac{2A}{i \hbar \omega } \tan \frac{t \omega}{2} 
\right]
\end{split}
\end{equation}

	Equation \eqref{eq:bufalos} is derived in \cite{omo:sprtopiq} derived in a different manner.  For the case $\wh{A} = \wh{I}$, it is a well known formula. (See, for instance, \cite{bay:dtq}.)


\section{Application: Bohr--Sommerfeld rules.} \label{sec:absr}

We explain in this section the application that caused our original interest in this problem.

  In \cite{col:bsrtao}, Colin de Verdi\'ere gives an algorithm which computes the Bohr--Sommerfeld quantization rules to all orders in $\hbar$ in the one dimensional case $N=1$.  His method is inspried by Voros \cite{vor:aheosqs} and a similar method had been previously used by Argyres \cite{arg:tbsqratwc}.  

   Let $\wh{H}$ be an operator with symbol $H \in C^{\infty}(T^*\rr)$.  Bohr--Sommerfeld quantization rules provide a way to asymptotically compute the spectrum of $\wh{H}$.
 Assume $H$ has a regular minimum at a point.  Under certain extra assumptions on the symbol $H$ (see \cite{col:bsrtao}), the eigenvalues of $\wh{H}$ are given by the solutions $E$ to 
\begin{equation}  \label{eq:action}
 2 \pi n \hbar \; = \; S(E) \; = \; \sum_{j=0}^{\infty} \hbar^j S_j(E)
\end{equation}
    with $n \in \zz$. To solve the previous equation, write $E=\sum \hbar^k E_k$ as a power series in $\hbar$ and substitute it into \eqref{eq:action} to obtain recursive expressions for each $E_k$.  This requires knowing the form of $S(E)$, called the \emph{semiclassical action,} for the hamiltonian $\wh{H}$.

  Let us consider for simplicity the case where the symbol $H$ does not depend on $\hbar$.  It is known that at lowest orders in $\hbar$:
\begin{equation*}
S_0(E) =  \int_{\gamma_E} p \; dx  \; , \quad \quad \quad S_1(E) = \pi
\end{equation*}
$\gamma_E$ denotes the level set $H^{-1}(E)$ around the minimum of $H$ and $(x,p)$ are the natural coordinates in $T^*\rr$.

      The main result in \cite{col:bsrtao} is the following:
\begin{equation*}
S_j(E) =  
\sum_{l=2}^{L(j)} \frac{(-1)^{l-1}}{(l-1)!} \left( \frac{d}{dE} \right)^{l-2}
\int_{\gamma_E} P_{j,l}(x,p) dt
\end{equation*}
where $t$ is the parametrization of $\gamma_E$ by the time evolution
\begin{equation*}
dx = H_p dt \; , \quad \quad \quad dp = - H_x dt
\end{equation*}
and $P_{j,l}$ are universal polynomials in the derivatives of $H$ defined by the symbol $R_a$ of the resolvent $\wh{R_a} = (a - \wh{A})^{-1}$:
\begin{equation*}
R_a = \frac{1}{a-H} \; + \; \sum_{j=1}^{\infty} \hbar^j
\sum_{l=2}^{L(j)} \frac{P_{j,l}(H)}{(a-H)^l}
\end{equation*}

Looking back at our formula for the symbol for the function of an operator \eqref{eq:primera} and using it for the function $f(y)=(a-y)^{-1}$, we see that $L_{j}=3j/2+1$  and the polynomials $P_{j,l}$ can actually be defined in terms of graphs:
\begin{equation*}
P_{j,l} (H) = \sum_{\substack{\textrm{reduced graphs } \Gamma \\ \textrm{with $l-1$ vertices}\\ \textrm{and $j$ edges} }}
\left( \frac{i}2 \right)^j \frac{c_\Gamma}{S_\Gamma} \lambda_\Gamma(H)
\end{equation*}
 which gives us the following equation for the eigenvalues $E$ of $\wh{H}$:
\begin{equation} \label{eq:goodBS}
\begin{split}
2 \pi (n & - 1/2)  \hbar \;\; = \;\;   S_0(E) \;  \\
&  + \; \sum_{\substack{\textrm{reduced graphs }\Gamma \\
\textrm{with }E_\Gamma > 0}}
\left( \frac{i \hbar}{2} \right)^{E_\Gamma} 
\frac{(-1)^{V_\Gamma}}{{V_\Gamma}!}
\left( \frac{d}{dE} \right)^{V_\Gamma -1}
\frac{c_\Gamma}{S_\Gamma}
\int_{\gamma_E} \lambda_\Gamma (H) \; dt
\end{split}
\end{equation}
 Here $E$ is an eigenvalue of $H$, whereas $E_\Gamma$ is the number of edges of a graph $\Gamma$.

$S_j(E)=0$ for $j>1$ and odd.  $S_2(E)$ is given by the contribution of 2 graphs, and $S_4(E)$ is given by the contribution of 15 graphs
     (see Appendix \ref{sec:c}).  However, there is a trick using Stokes' theorem that allow us to express the contribution of certain graphs in this expression in terms of others.  (This trick is used in \cite{car:mspattbsa} and in \cite{col:bsrtao} for $S_2$, although without the diagrammatic notation.)
 As a consequence, $S_2(E)$ can be written in terms of 1 graph and $S_{4}(E)$ can be written in terms of 5 graphs (those where every vertex has at least two edges):
\begin{equation} \label{eq:BS4}
\begin{split}
2 \pi &(n  -1/2) \hbar \; \; = \; \; S_0(E) \; \; \\
&  - \frac{\hbar^2}{4} \; \Bigg[ \; 
     \frac{1}{2! \cdot 6} \; \frac{d}{dE} \int_{\gamma_E} \! \! \!
     \xymatrix{(H \ar[r] \ar@<-.5ex>[r] & H)}
     dt \; \Bigg] \\
& + \frac{\hbar^4}{16} \Bigg[ \;    
     \frac{1}{2! \cdot 120} \; \frac{d}{dE} \; \int_{\gamma_E} \! \!
\xymatrix{(H \ar[r] \ar@<.5ex>[r] \ar@<-.5ex>[r] \ar@<-1ex>[r] & H)} dt
   \; + \; 
     \frac{1}{4! \cdot 12} \left( \frac{d}{dE} \right)^3 \int_{\gamma_E} \! \!
\xymatrix{(H \ar[r] \ar@<-.5ex>[r] & H)^2} dt  \\
&  \quad \quad - \; 
     \frac{1}{3! \cdot 15} \left( \frac{d}{dE} \right)^2 \! \int_{\gamma_E} \! 
\xymatrix{H \ar[r] \ar@<-.5ex>[r] \ar[rd] & H  \ar[d] \\ & H } dt
  \; + \;
     \frac{1}{4! \cdot 15} \left( \frac{d}{dE} \right)^3 \int_{\gamma_E} \! 
\xymatrix{H \ar[r] & H \ar[d] \\ H \ar[u] & H \ar[l]} dt  \\
& \quad \quad - \; 
     \frac{1}{3!\cdot 12} \left( \frac{d}{dE} \right)^2 \int_{\gamma_E} \! \!
\xymatrix{(H \ar[r] \ar@<-.5ex>[r] & H \ar[r] \ar@<-.5ex>[r] & H)} dt
\; \Bigg] \\
& + O(\hbar^6)
\end{split}
\end{equation}

   All the integrands in the previous expressions are long polynomials in the derivatives of $H$.  The expression would be hard to obtain without the diagrammatic notation.  Given a concrete hamiltonian $H$ we could easily program a computer to write all the terms in \eqref{eq:goodBS} for that specific operator at higher orders in $\hbar$.

   For the case of a hamiltonian of the form kinetic plus potential energy
$H(x,p) = p^2/2m + V(x)$ the contribution of many graphs vanishes, and \eqref{eq:BS4} becomes:
\begin{equation*}
\begin{split}
2 \pi &(n-1/2) \hbar \; \; = \; \; S_0(E) \; \;
- \frac{\hbar^2}{m} \frac{1}{24} \frac{d}{dE} \int_{\gamma_E} V''(x) dt \\
& +  \frac{\hbar^4}{m^2} \frac{1}{2^7 3^2} \left[ \;
\frac{7}{5} \left( \frac{d}{dE} \right)^3 \! \int_{\gamma_E} \! \left[ V''(x) \right]^2 dt \; 
- \left( \frac{d}{dE} \right)^2 \! \int_{\gamma_E} \! V^{(4)}(x) dt \;
\right] \;  + O \left(\frac{\hbar^6}{m^3} \right)
\end{split}
\end{equation*}
    
The disadvantage of this method is that it does not generalize to the multidimensional case $N > 1$.  Cargo et al
\cite{car:mspattbsa} approached this problem in a totally different way to obtain a result valid in all dimensions.  The symbol of a function of an operator plays a role in their derivation, too.


\appendix


\section{A lemma for calculations with graphs.}  \label{sec:slfcwg} \label{sec:alfcwg}

The diagramatic notation introduced in \textsection \ref{sec:gra} makes 
equations and derivations easier to write.  The following lemma 
is needed in some of those derivations:

\begin{lem} \label{lem:diag1}
Let $\Gamma$ be a labeled graph with $V$ vertices and $E$ edges.  Let $D, C_1, \ldots , C_V$ be symbols.  Then
	\begin{equation*}
	\xymatrix{ \lambda_\Gamma (C_1, \ldots C_V) 
\ar@<1.5ex>[r] \ar@<-1.5ex>[r]^(.7){\ka} & D}  \; = \; \;\sum_{\Gamma'} \; 
\lambda_{\Gamma'} (C_1, \ldots, C_V, D)
\end{equation*}
	where the sum is taken over all labeled graphs $\Gamma'$ with $V'=V+1$ vertices and $E'=E+k$ edges, 
which are constructed by putting together
\begin{itemize}
\item[-] the labeled graph $\Gamma$ (conserving its labels),
\item[-] an extra vertex labeled by $V+1$, 
\item[-] $k$ extra arrows (labeled by $E+1,\ldots,E+k$) starting from the vertices of $\Gamma$ and ending at the vertex $V+1$.
\end{itemize}
 \end{lem}
        \begin{proof}  Write down the definition of both sides in terms of 
$J^{\mu \nu}$ and check that they are equal.
\end{proof}

For instance:
\begin{equation*}
\begin{split}
& \xymatrix{(C \ar[r]_{1} & D) \ar@<1ex>[r] \ar[r] & E} \quad  \\ 
\quad & \; = \;
\xymatrix{C \ar[r]_{1} \ar@/^/[rr]^{2} \ar@/^/@<.5ex>[rr]_(.7){3} & D & E}
\quad + \quad
\xymatrix{C \ar[r]_{1} & D \ar@<1ex>[r]^{2} \ar[r]_{3} & E}  \quad  \\
\quad & \quad  \quad + \;
\xymatrix{C \ar[r]_{1} \ar@/^/[rr]^{2} & D  \ar[r]_{3} & E} \quad + \quad
\xymatrix{C \ar[r]_{1} \ar@/^/[rr]^{3} & D  \ar[r]_{2} & E} \quad  \\
\quad & = \;
\xymatrix{C \ar[r] \ar@/^/[rr] \ar@/^/@<.5ex>[rr] & D & E}
\quad + \quad
\xymatrix{C \ar[r] & D \ar@<1ex>[r] \ar[r] & E}  \quad  \\
\quad & \quad \quad + \; 2 \;
\xymatrix{C \ar[r] \ar@/^/[rr] & D  \ar[r] & E} 
\end{split}
\end{equation*}

\section{An alternative derivation of \eqref{eq:function}.}  \label{sec:aado}

As we mentioned in the introduction and in \textsection \ref{sec:absr},
 the existence of a universal formula like 
\eqref{eq:primera} was used by Voros \cite{vor:aheosqs} and 
Colin de Verdi\`ere 
\cite{col:bsrtao} (and, indirectly, by Argyres \cite{arg:tbsqratwc}).
They start by writing a smooth function of an operator in terms 
of the resolvent.
Let $a \in \cc$ and define the resolvent operator
$\wh{R_a} : = (a - \wh{A})^{-1}$ with symbol $R_a$.  
Then we use Helffer--Sj\"ostrand's formula \cite{hel:edsacmeedh} as a spectral theorem:
\begin{equation*} 
\wh{B} \; = \; f(\wh{A}) \; = \; -\frac{1}{\pi} \; 
\int_{\cc} \wh{R_z} \;
  \partial_{\overline{z}} \widetilde{f} (z) \; dx dy
\end{equation*}
      Here $z=x+iy$, $\widetilde{f}$ is an almost analytic extension of $f$,
and $\partial_{\overline{z}}=\partial_x + i \partial_y$.  This allows us 
to write for the symbol:
\begin{equation} \label{eq:HS}
B \; = \; -\frac{1}{\pi} \; 
\int_{\cc} R_z \;
  \partial_{\overline{z}} \widetilde{f} (z) \; dx dy
\end{equation}
    
Hence, finding the symbol of $f(\wh{A})$ reduces to finding the 
symbol of $\widehat{R_a}$.\footnote{We could also use Cauchy's integration formula as a spectral theorem:
\begin{equation*}
\begin{split}
\wh{B} \; & = \; f(\wh{A}) \; = \; \int_{\gamma} \frac{da}{2\pi i} \; f(a)
\wh{R_a} \\
B \; & =  \; \int_{\gamma} \frac{da}{2\pi i} \; f(a) R_a 
\end{split}
\end{equation*}
which will lead to the same results, but it is only valid for analytic functions $f$.  $\gamma$ is a path around the spectrum of $\wh{A}$.}
  In order to do so, we may write
     \begin{equation*}
R_a \; = \; \sum_{k=0}^{\infty} R_{a (k)} \hbar^k
\end{equation*} 
and, since $\wh{R_a}(a-\wh{A})=1$, substitute it into
\begin{equation*}
R_a \star (a-A) \; = \; (a-A) \star R_a = 1
\end{equation*}
to obtain recursively the value of each $R_{a(k)}$.
	Although simple, this method quickly proves intractable.  The 
calculations at order $4$ are already too complex and we will not find the
pattern that \eqref{eq:primera} shows. 

	But we can also use this idea to prove \eqref{eq:function}
in a different way.  First, we prove it for the function 
$f(y) = (a-y)^{-1}$:

\begin{note}[Claim.] The function
\begin{equation} \label{eq:redres}
h_a(A) = \sum_{\substack{\text{reduced}\\
        \text{labeled} \\ \text{graphs } \Gamma}}   
        \frac{1}{E!} \left(\frac{i \hbar}{2}\right)^E
        \frac{\lambda_\Gamma(A)}{(a-A)^{V+1}}
\end{equation}
satisfies $h_a(A) \star (a-A) = (a-A) \star h_a(A) = 1$.
\end{note}
The proof is a long combinatorial exercise on calculations with graphs.
Therefore $h_a(A) = R_a$.

Second, we can substitute \eqref{eq:redres} into 
\eqref{eq:HS} to  
obtain again \eqref{eq:function}.


\section{Calculation of $c_\Gamma$ and $S_\Gamma$.}
 \label{sec:coc} \label{sec:c} \label{sec:cocas}

	Calculating $c_\Gamma$ is a combinatorial problem.
The following five facts give quick, recursive rules for it:

\begin{note}[Facts:]
\begin{enumerate}
	\item \label{it:1}
$c_\isol = 1$
	\item \label{it:odd} $c_\Gamma = 0$ if $E$ is odd.

	This is due to a symmetry property.  If we denote by $\Gamma$ a 
labeled
graph with a numbering of the vertices by $1,2,\cdots V$ and by 
$i(\Gamma)$ the relabeling of the vertices by the permutation 
$\binom{1 2 \cdots V}{V \cdots 2 1}$,
 then $$\lambda_{i(\Gamma)} = (-1)^E \lambda_\Gamma$$
	And when we sum over reorderings of the vertices:
$$c_{\Gamma} = (-1)^E c_{\Gamma}$$
	\item \label{it:dis}
If $\Gamma_i$ has $V_i$ vertices, $\Gamma_i \neq \Gamma_j$ for $i \neq j$, and $\Gamma = \Gamma_1 + 
\cdots \Gamma_n$ is the topological sum with $V=V_1 + \cdots V_n$ 
vertices, then
\begin{equation} \label{eq:cdis}
	c_\Gamma = \frac{V!}{V_1! \cdots V_n!} c_{\Gamma_1} \cdots 
c_{\Gamma_n}
\end{equation}
	\item $c_\Gamma$ does not change if we erase two edges with the 
same endpoints.

	For instance $$c_{1 \ra 2 \rra 3 \ra 4} = 
c_{1 \ra 2 \quad 3 \ra 4}$$
	\item \label{it:p}
If $p$ is a vertex in $\Gamma$, denote by $\Gamma - p$ the 
same $\Gamma$ with the vertex $p$ and every edge starting or ending at 
$p$ erased.  
		 For instance:
\begin{equation*}
\begin{split}
         \Gamma \quad = & \quad
        \xymatrix{ \bullet \ar[r] & p \ar[r] & \bullet
        \ar[r] \ar@<-.5ex>[r] & \bullet} \\
\Gamma - p \quad = & \quad
        \xymatrix{ \bullet  &   & \bullet
        \ar[r] \ar@<-.5ex>[r] & \bullet}
\end{split}
\end{equation*}
Then:
\begin{equation*}
	c_\Gamma = \sum_{\text{vertices $p$ in $\Gamma$}} (-1)^{\text{\# 
of arrows \emph{starting} at p}} \;  c_{\Gamma - p}
\end{equation*}
\begin{proof}
	Count the possible reorderings of the vertices by choosing first 
which vertex has label $V$.
\end{proof}
\end{enumerate}
\end{note}

In particular, putting together facts \ref{it:odd} and \ref{it:dis} we
get:
\begin{lem} \label{lem:odd}
If  $\Gamma$ has a connected component with an odd number of edges, then $c_\Gamma = 0$.
\end{lem}

These are the reduced, connected graphs with 
2 and 4 edges, and for each of them the value of $c_\Gamma$ and 
$S_\Gamma$:

\begin{equation*}
\begin{array}{|c|c|c|c|c|}
\hline
	\Gamma & E & V  & S_\Gamma  & c_\Gamma \\
\hline \hline
	\xymatrix{\bullet \ar[r] \ar@<-.5ex>[r]& \bullet} & 2&2& 4 & 2 \\
\hline
	\xymatrix{\bullet \ar[r] & \bullet \ar[r] & \bullet} &2&3&2&-2 \\
\hline
	\xymatrix{\bullet \ar[r] \ar@<-.5ex>[r] \ar@<-1ex>[r] 
\ar@<.5ex>[r] & \bullet} & 4 &2& 48 & 2 \\
\hline
	\xymatrix{\bullet \ar[r] \ar@<-.5ex>[r] \ar@<.5ex>[r] & \bullet
\ar[r] & \bullet} & 4 & 3&6 & -2 \\
	\xymatrix{\bullet \ar[r] \ar@<-.5ex>[r] & \bullet \ar[r] 
\ar@<-.5ex>[r] & \bullet} &&& 8 & 6 \\
	\xymatrix{\bullet \ar[r] \ar@<-.5ex>[r] \ar[rd] & \bullet
\ar[d] \\ & \bullet} &&& 4 & 2 \\
\hline
	\xymatrix{\bullet \ar[r] & \bullet \ar[d] \\ \bullet \ar[u] &
\bullet \ar[l]} &4&4& 8 & 8 \\
	\xymatrix{\bullet \ar[r] & \bullet \ar[r] \ar[d] & \bullet \\
& \bullet \ar[ul]} &&& 2 & 0 \\ 
	\xymatrix{\bullet \ar[r] \ar@<-.5ex>[r] & \bullet \ar[r] \ar[d] 
& \bullet \\ & \bullet} &&& 4 & 8 \\
	\xymatrix{\bullet \ar[r] \ar@<-.5ex>[r] & \bullet \ar[r]
& \bullet \ar[r] & \bullet} &&& 2 & -8 \\
	\xymatrix{\bullet \ar[r] & \bullet \ar[r] \ar@<-.5ex>[r] & \bullet 
\ar[r] & \bullet} &&& 4 & 0 \\
\hline
	\xymatrix{\bullet \ar[r] & \bullet \ar[r] & \bullet \ar[r] &
\bullet \ar[r] & \bullet} & 4 & 5 & 2 & 16 \\
	\xymatrix{\bullet \ar[r] & \bullet \ar[r] & \bullet \ar[r] 
\ar[d] & \bullet \\ && \bullet} &&& 2 & -8 \\
	\xymatrix{\bullet \ar[r] & \bullet \ar[r] \ar[d] \ar[dr] &
\bullet \\ & \bullet & \bullet} &&& 24 & -24 \\
\hline
\end{array}
\end{equation*}
  

\section{Symbol of a function of an operator at order 4 in $\hbar$.}     
  \label{sec:soafoaoao4ih}

  Using any of the equations derived in \textsection \ref{sec:cal}, we write down the explicit form of all the terms  of the symbol of a function of an operator up to order $4$ in $\hbar$.  The data in the table in Appendix \ref{sec:cocas} are needed.

\begin{equation*}
\begin{split}
	\widehat{B} \; = & \; f(\widehat{A}) \\
	B \; = & \; f(A)  \\
	& - \frac{\hbar^2}{4} \; \Bigg[ \;
\frac{\xymatrix{A \ar[r] \ar@<-.5ex>[r]& A}}{2} \; 
\frac{f''(A)}{2!} \quad + 
	\quad \xymatrix{A \ar[r] & A  & A \ar[l]} \; 
\frac{f'''(A)}{3!} \;
\Bigg] \\
	& + \; \frac{\hbar^4}{16} \Bigg[ \;
	\frac{\xymatrix{A \ar[r] \ar@<-.5ex>[r] \ar@<-1ex>[r]
\ar@<.5ex>[r] & A}}{24} \; \frac{f''(A)}{2!} \\
	& \quad \quad \quad + \Bigg(
	 \frac{\xymatrix{ A \ar[r] \ar@<-.5ex>[r] \ar@<.5ex>[r] &
A & A \ar[l]}}{3} \quad + \quad \frac{1}{2} \;
	\xymatrix{A \ar[r] \ar@<-.5ex>[r] \ar[rd] & A \ar[d] \\ & A} \\
& \quad \quad \quad \quad \quad
	+ \; \frac{3}{4} \; 
	\xymatrix{A \ar[r] \ar@<-.5ex>[r] & A \ar[r] \ar@<-.5ex>[r] &A}
\; \quad \quad \quad \quad \quad \quad \quad \quad 
	\Bigg) \frac{f'''(A)}{3!} \\
	& \quad \quad \quad + \Bigg( \;
	\frac{3}{4} \; (\xymatrix{A \ar[r] \ar@<-.5ex>[r]& A})^2 \quad
+ \quad \xymatrix{A \ar[r] & A \ar[d] \\ A \ar[u] & A \ar[l]} \\
& \quad \quad \quad \quad \quad
	+ \; 4 \; \xymatrix{A \ar[r] \ar@<-.5ex>[r] & A \ar[r] & A
& A \ar[l]} \\
& \quad \quad \quad \quad \quad
	+ \; 2 \; \xymatrix{A \ar[r] \ar@<-.5ex>[r] & A
\ar[r] \ar[d] &A \\ &A} \;
\quad \quad \quad \quad \quad 	
	\Bigg) \; \frac{f^{(4)}(A)}{4!} \\
	& \quad \quad \quad + \Bigg( 
	8 \;\xymatrix{ A \ar[r] & A \ar[r] & A \ar[d] \\ & A & A \ar[l]}
	\quad + \quad  
	\xymatrix{A \ar[r] &A & A \ar[l] \\ &A \ar[u] &A \ar[ul]} \\
& \quad \quad \quad \quad \quad
	+ \; 5 \; (\xymatrix{A \ar[r] \ar@<-.5ex>[r]& A}) \;
	(\xymatrix{A \ar[r] & A & A \ar[l] }) \\
& \quad \quad \quad \quad \quad
	+ \; 4 \; 
	\xymatrix{ A & A \ar[l] \ar[r] & A \ar[d] \ar[r] &A \\ &&A} \;
\quad \quad \quad \quad \quad 
	\Bigg)\; \frac{f^{(5)}(A)}{5!} \\
	& \quad \quad \quad + \; 10 \; ( \xymatrix{A \ar[r] & A  & 
A \ar[l]})^2 \;
	 \quad \frac{f^{(6)}(A)}{6!} \; \Bigg] \\
	& + O(\hbar^6)
\end{split}
\end{equation*}


\bibliographystyle{amsplain}

\bibliography{bibliography}

\end{document}